\documentclass[12pt]{article}
\usepackage{mathbbold}
\usepackage{float}
\usepackage{amsfonts}
\usepackage{
mathrsfs,amsmath}
\usepackage{amssymb}
\usepackage{epsfig}
\usepackage{subfigure}
\usepackage{fancybox}
\usepackage{graphicx}
\usepackage{enumerate}
\newcommand{\diag}{{\rm diag}}
\newcommand{\lbar}{\overline}
\newcommand{\wdh}{\widehat}

\def\op{{\cal L}}
\def\l{\left|}
\def\r{\right|}

\newcommand{\wdt}{\widetilde}
\newcommand{\e}{\varepsilon}
\newcommand{\rr}{{\Bbb R}}
\newcommand{\M}{{\cal M}}
\newcommand{\cd}{(\cdot)}
\newcommand{\nd}{\noindent}

\def\para#1{\vskip .4\baselineskip\noindent{\bf #1}}
\def\qed{\strut\hfill $\Box$}

\newtheorem{thm}{Theorem}[section]

\newtheorem{lem}[thm]{Lemma}

\newtheorem{exm}[thm]{Example}

\newcommand{\thmref}[1]{Theorem~{\rm \ref{#1}}}
\newcommand{\lemref}[1]{Lemma~{\rm \ref{#1}}}

\newcommand{\exmref}[1]{Example~{\rm \ref{#1}}}

\def\al{\alpha}

\newcommand{\beq}[1]{\begin{equation} \label{#1}}
\newcommand{\eeq}{\end{equation}}
\newcommand{\bed}{\begin{displaymath}}
\newcommand{\eed}{\end{displaymath}}
\newcommand{\bea}{\bed\begin{array}{rl}}
\newcommand{\eea}{\end{array}\eed}
\newcommand{\ad}{&\!\!\!\disp}
\newcommand{\aad}{&\disp}
\newcommand{\barray}{\begin{array}{ll}}
\newcommand{\earray}{\end{array}}
\def\({\left(}
\def\){\right)}
\def\disp{\displaystyle}

\numberwithin{equation}{section}

\begin{document}

\title{A Mean-Variance Control Framework
for Platoon Control Problems: Weak Convergence Results
 and Applications on Reduction of Complexity\thanks{This research was
supported in part by the National Science Foundation under
CNS-1136007.}}
\author{Zhixin Yang,\thanks{Department of Mathematics,
Wayne State University, Detroit, Michigan 48202, 
zhixin.yang@wayne.edu.} \and G. Yin,\thanks{Department of
Mathematics, Wayne State University, Detroit, Michigan 48202, 
gyin@math.wayne.edu.} \and Le Yi Wang,\thanks{Department of
Electrical and Computer Engineering, Wayne State University,
Detroit, MI 48202, 
lywang@wayne.edu.} \and Hongwei
Zhang\thanks{Department of Computer Science, Wayne State University,
Detroit, MI 48202,
hongwei@wayne.edu.}}
\date{Dedicated to Wing Shing Wong on the occasion of his 60th Birthday}
\maketitle

\begin{abstract}

This paper introduces  a new approach of treating platoon systems
using mean-variance control formulation. The underlying system is a
controlled switching
diffusion in which the random switching process is a continuous-time
Markov chain. This switching process is used to represent random
environment and other random factors that cannot be given by stochastic
differential equations driven by a Brownian motion.
The state space of the
Markov chain is large in our setup, which renders practically infeasible a
straightforward implementation
of the mean-variance control strategy obtained in the literature. By partitioning
 the states of the Markov chain
 into sub-groups (or clusters) and then aggregating the states of each cluster as a
super state, we are able to obtain a limit system of much reduced complexity. The justification of the
limit system is rigorously supported by establishing certain weak convergence
results.

\medskip
\nd{\bf Key Words.} platoon control, mean-variance control, two-time-scale model, weak convergence, reduction of
complexity.
\end{abstract}

\newpage

\setlength{\baselineskip}{0.25in}
\section{Introduction}

Highway vehicle control is a critical task in developing
intelligent transportation systems.
Platoon formation has been identified
as one promising strategy for enhanced safety,
improved highway utility,
increased fuel economy, and reduced emission
toward autonomous or semi-autonomous vehicle control.
The goal of longitudinal platoon control is
to ensure that all the vehicles move in
the same lane at the same speed with desired inter-vehicle distances.

Platoon control has been studied in the contexts of intelligent
highway control and automated
highway systems for many years with numerous methodologies and demonstration
systems \cite{HMS,RTLZ}. Many control methodologies have been applied,
including PID controllers, state feedback, adaptive control, state observers,
among others, with safety, string stability, and team coordination
as the most common objectives \cite{CH1,LP,SH}.

A platoon may be viewed as a networked system consisting of
individual subsystems whose operations and resource
consumptions must be carefully coordinated to achieve
desired performance for  the entire platoon.
Platoon control and performance optimization bear
certain similarity to portfolio management
in mathematical finance in which optimal distribution
of available resources to different subsystems
(stocks or mutual funds) can lead to increased return and reduced risk.
By borrowing the basic idea of mean-variance
control from mathematical finance,
we study its potential extension and applications to
platoon control systems.

The origin of the mean-variance optimization problem can be traced
back to the Nobel-prize-winning work of Markowitz \cite{Mark52}. The
salient feature of the model is that, in the context of finance, it
enables an investor to seek highest return after specifying the
acceptable risk level quantified by the variance of the return. The
mean-variance approach has become the foundation of modern finance
theory and has inspired numerous extensions and applications. Using
the stochastic linear-quadratic (LQ) control framework, Zhou and Li
\cite{ZhouL} studied the mean-variance problem for a continuous-time
model. Note that the problem becomes fundamentally different from
the traditional LQ problem studied in literature. In the classical
time-honored LQ theory, the matrix related to the control (known as
control weight) needs to be positive definite. In the mean variance
setup for linear systems, the control weight is non-positive
definite. In our previous work
\cite{ZhouY}, the mean-variance problems for switching
diffusion models were treated and a number of results including
optimal portfolio selection,
efficient frontier, and mutual fund theory were discovered.

In this study, we identify the following
three scenarios in platoon control problems
in which resource allocation and risk management
lead naturally to mean-variance formulations.
\begin{enumerate}
\item
 Consider the longitudinal
 inter-vehicle distance control. To increase
 highway utility, it is desirable
 to reduce the total length of a platoon,
 which intends to reduce inter-vehicle distances.
 This strategy, however, will
 increase the risk of collision in the presence of
 vehicle traffic uncertainties.
 This tradeoff amounts to maximizing benefits at a tolerable risk. This may be compared
 to financial portfolio
 management problems in which one wants to maximize profit return but control
 the risk too. Consequently, the basic idea of mean-variance (MV) control becomes useful.
 The MV approach has never been applied to platoon control. It offers several distinct advantages:
 1) Unlike heuristic methods such as neural network optimization and genetic algorithms,
  the MV method is simple but rigorous; 2) the MV method is computationally efficient;
   3) the form of the solution (i.e., efficient frontier) is readily applicable
   to assessing risks in platoon formation, hence is practically appealing.

 \item Consider communication resource
 allocation of bandwidths for vehicle to vehicle (V2V)
 communications. For a given maximum throughput
 of a platoon communication system,
 the communication system operator
 must find a way to assign this resource
 to different vehicle-to-vehicle channels.
 Each channel's bandwidth usage is the state
 of the subsystem.
 Their summation is a random process
 and is desired to approach the maximum throughput
 (the desired mean at the terminal time)
 with small variations.
 Consequently, it becomes a mean-variance control problem.

\item
 We may view platoon fuel consumption
 (or similarly, total emission) in the MV setting.
 The platoon fuel consumption is
 the summation of vehicle fuel consumptions.
 Due to variations in vehicle sizes and speeds,
 each vehicle's fuel consumption is a controlled
 random process.
 Tradeoff between a platoon's team acceleration/maneuver
 capability and fuel consumption can be summarized as
 a desired platoon fuel consumption rate.
 Assigning allowable fuel consumption rates  to different vehicles result
 in coordination of vehicle operations modeled
 by subsystem fuel rate dynamics.
 To control the platoon fuel consumption
 rate to be close to the designated value,
 one may formulate this as a mean-variance control problem.
\end{enumerate}

Due to vehicle mobility and network resource fluctuations,
platoon network topologies
or random environment may vary dynamically.
To capture this common feature in platoon control,
we model the network topology or the environment as
a continuous-time Markov chain. The resulting system becomes one having
regime switching.  We assume that the
 Markov chain
 has a large state
space in order to deal with complex systems.
To treat the platoon problems, we could in
principle apply the results in \cite{ZhouY}.
Nevertheless, the large state space of the Markov chain
renders a
straightforward implementation
of the mean-variance control strategy obtained in \cite{ZhouY}
practically infeasible. The computational complexity becomes a major
concern.
Inspired by the idea in the work  \cite{SimonA},
to exploit the hierarchical structure of the underlying systems,
and to fully utilize the near decomposability
\cite{Courtois,SethiZ94}
by means of considering fast and slow switching modes,
the work \cite{ZhangY} treated near-optimal control problems of LQG
with regime switching. However, only positive definite control
weights were allowed there under the usual quadratic control criteria.
In our current setup, the control weights
are indefinite, so the main assumption in \cite{ZhangY}  does not hold.
Physically, only part of the network topology or random environment
will change at a time such as the addition/departure
of a vehicle, or the loss/addition of a
communication link. Some parts of the chain vary
rapidly (e.g., when vehicles pass some bridges
that block signal  transmissions) and
others change slowly (e.g., when vehicles are
moving smoothly in open space). The fast and slow variations are in high
contrast, resulting in a two-time-scale formulation.
This paper, together with its companion paper \cite{YYWZ}, sets up a new formulation
towards resolving platoon coordination and optimization issues with
reduced computational complexity.

This two-time-scale scenario provides an opportunity
to reduce computational complexity for the Markov chain.
The main idea is a decomposition of the large space
into sub-clusters and aggregation of states in each sub-cluster. That is,
we partition the state space of the Markov chain into subspaces (or sub-groups or
sub-clusters). Then,
in each of the sub-clusters, we aggregate all the states into
one super state. Thus the total number of discrete states is substantially
reduced. In the companion paper \cite{YYWZ}, we obtained near-optimal controls
by designing controls using a ``limit system'' and such constructed controls are
nearly optimal.
This paper focuses on justifying the limit system by means of weak convergence
methods. The weak convergence result is proved using a martingale
problem formulation.
As a first step in this direction,
this paper presents the key mathematical framework and core results.
Their usage can be expanded by considering further
system details in practical applications with concrete
model structures, sizes, resource definitions,
and physical limitations. These will
be investigated and reported elsewhere.

The rest of the paper is arranged as follows.
Section \ref{sec:form} formulates the
 two-time-scale platoon problems. Section \ref{sec:mv}
proceeds with the study of the underlying mean-variance problem.
 Section \ref{sec:weak}
derives weak convergence of underlying systems using martingale
problem formulation, which rigorously justifies the use of the limit
system. The limit system has a substantially fewer number of
discrete states, resulting in much reduced complexity. Section
\ref{sec:near} recalls the near-optimal controls obtained in
\cite{YYWZ} and  presents numerical experiments that further
illustrate the near optimality. Finally, Section \ref{sec:fur}
concludes the paper with some final remarks.

\section{Problem Formulation}\label{sec:form}
We work with a complete probability space $(\Omega,
\mathcal{F}, P)$. Suppose that $\al(t)$ is
  continuous-time  Markov
chain with state space $\mathcal{M}=\{1,2,\ldots,m\}$, that
$w(t)=(w_1(t),w_2(t),\ldots,w_d(t))'$ is a standard $d$-dimensional
Brownian motion, where $a'$ denotes the transpose of $a\in
\mathbb{R}^{l_1 \times l_2}$ with $l_i \ge 1$, and
$\mathbb{R}^{r\times 1}$ is simply written as $\mathbb{R}^r$ in what
follows. Suppose that  $w(t)$ and the Markov chain are independent
of each other.
In \cite{ZhouY}, a mean-variance portfolio selection problem in
which the environment is randomly varying and modeled by  a
regime-switching system was treated. In this paper, we continue to use
the same setup as in \cite{ZhouY}. In addition to the finance
applications, we are particularly interested in platoon control
problems. Mathematically, the new feature considered here is that
the state space of the discrete event process $\al\cd$ is large.
Obtaining the optimal strategy in such a large-scale system involves
high computational complexity,
optimal control a difficult task. To reduce the computational
complexity, we note that in the Markov chain, some groups of states
vary rapidly whereas others change slowly. Using the distinct
transition rates, we decompose the state space $\M$ into subspaces
$\M = \cup^l_{i=1} \M_i$ such that within each $\M_i$, the
transitions happen frequently and among different clusters the
transitions are relatively infrequent. To reflect the different
transition rates, we let $\al(t)=\al^\e(t)$ where $\e>0$ is a small
parameter so that the generator of the Markov chain is given by
\beq{mar}%
Q^\e=\frac{\wdt Q}{\e}+\wdh{Q}. \eeq
We consider a network that contains  $d_1+1$
nodes. The flow of one of the nodes is given by
stochastic ODE 
\beq{0}\barray dx_0^\e(t)\ad=r(t,\al^\e(t))x_0^\e(t)dt\\
x_0^\e(0)\ad=x_0, \al^\e(0)=\al.\earray\eeq 
 where for each $j\in \M$, $r(t,j)\ge 0$ is 
the growth rate corresponding to $x_0^\e(t)$. The
flows of the other $d_1$ nodes follow geometric Brownian motion models
under random environments (random switching) \beq{2}\barray
dx_i^\e(t)\ad=x_i^\e(t)r_i(t,\al^\e(t))dt+x_i^\e(t)\sigma_i(t,\al^\e(t))dw(t)\\
x_i^\e(0)\ad=x_i,\al_i^\e(0)=\al. \text{ for }
i=1,2,\ldots,d_1.\earray \eeq where
\bea
\ad\sigma_i(t,\al^\e(t))=(\sigma_{i1}(t,\al^\e(t)),\sigma_{i2}(t,\al^\e(t)),
\ldots,\sigma_{id}(t,\al^\e(t)))\in\rr^{1\times
d},\\
\ad w(t) =(w_1(t),w_2(t),\ldots,w_d(t))'\in\rr^{d\times 1},\eea
and $r_i(t,j)\in \rr$ (with $i=1,\ldots, d_1$ and 
$j\in \M$) is the drift rate for the flow of
the ith node. We can represent the total flows of the whole
network system as  a 1-dimensional variable $x^\e(t)$ for which
we need to decide the proportion $n_i(t)$ of flow $x_i^\e(t)$ to put
on node $i$, i.e., \bea \disp
x^\e(t)=\sum^{d_1}_{i=0}n_i(t)x_i^\e(t).\eea By assuming that the
interaction among these $d_1+1$ nodes occurs continuously, we
have
\beq{x}%
\barray%
dx^\e(t)\ad=\sum^{d_1}_{i=0}n_i(t)dx_i^\e(t)\\
\ad=[r(t,\al^\e(t))x^\e(t)+B(t,\al^\e(t))u(t)]dt+u'(t)\sigma(t,\al^\e(t))dw(t)\\
x^\e(0)\ad=x=\sum^{d_1}_{i=1}n_i(0)x_i , \ \al^\e(0)=\al, \text{ for }
0 \le t \le T,\earray \eeq where
\bea \ad  B(t,\al^\e(t))=(r_1(t,\al^\e(t))
-r(t,\al^\e(t)),r_2(t,\al^\e(t))-r(t,\al^\e(t)),
\ldots,\\
\aad \quad \hfill r_{d_1}(t,\al^\e(t))-r(t,\al^\e(t))),\\
\ad \sigma(t,\al^\e(t)) =
(  \sigma _1 (t,\al^\e(t)),\ldots, \sigma _{d_1}
  (t,\al^\e(t)))'\in\rr^{{d_1}\times d},\\
\ad   u(t)=(u_1(t),\ldots, u_{d_1}(t))'\in\rr^{d_1\times 1},\eea
 and
$u_i(t)=n_i(t)x_i(t)$
 is the total amount of flow
for node $i$ at time $t$ for $i=1,2,\ldots,d_1$.
We assume throughout this paper that all the functions $r(t,i)$,
$B(t,i)$, and $\sigma(t,i)$ are measurable and uniformly bounded in
$t$.
We also assume the non-degeneracy condition is satisfied, i.e.,
there is a $\delta >0$ such that $a(t,i)=\sigma(t, i)\sigma'(t,i)\ge
\delta I$ for any $t\in[0,T]$ and $i\in \M$. We denote by
$L^2_{\mathcal{F}}(0,T;\mathbb{R}^{l_0})$ the set of all
$\mathbb{R}^{l_0}$-valued, measurable stochastic processes $f(t)$
adapted to $\{\mathcal{F}_t\}_{t\ge 0}$ such that
$E\int^T_0|f(t)|^2dt<+\infty$.

Let $\mathcal{U}$ be the set of controls which is a
compact set in $\mathbb{R}^{d_1\times 1}$. The $u(\cdot)$ is said to
be admissible if  $u(\cdot)\in
L^2_{\mathcal{F}}(0,T;\mathbb{R}^{d_1})$ and the equation \eqref{x} has
a unique solution $x^\e(\cdot)$ corresponding to $u(\cdot)$. In this
case, we call $(x^\e(\cdot),u(\cdot))$ an admissible
(total flow, flow distribution) pair.
Our objective is to find an admissible control $u(\cdot)$ among all
the admissible controls given that the expected terminal flow value of
the whole system is $Ex^\e(T)=z$ for some given $z\in \mathbb{R}$ so
that the risk measured by the variance at the terminal of the flow is
minimized. Specifically, we have the following performance measure \beq{obj}
\barray \ad \min
\left\{J(x,\al,u(\cdot))=E[x^\e(T)-z]^2\right\}\\
\ad \text{ subject to } Ex^\e(T)=z. \earray \eeq Note that
in this case, the objective function does not involve
control $u$. Thus, the
LQG problem
is one with zero control weight hence the problem
becomes one with indefinite control weights.

\section{Feasibility and Optimal Controls}\label{sec:mv}
To begin, we present the following
lemma, whose proof can be found in \cite[Theorem 3.3]{ZhouY}.

\begin{lem}
The mean variance problem \eqref{obj} is feasible for every $z\in
\rr$ if and only if
$$E\left\{\int^T_0|B(t,\al^\e(t))|^2dt\right\}>0.$$
\end{lem}

To study optimality and to handle the
constraint in \eqref{obj}, we apply the Lagrange
multiplier technique and get unconstrained problem (see,
e.g.,\cite{ZhouY}) with multiplier $\lambda\in \mathbb{R}$:
\beq{obj2}%
\barray%
 \ad \min\left\{
J(x,\al,u(\cdot),\lambda)=E[x^\e(T)+\lambda-z]^2-\lambda^2\right\}\\
\ad \text{ subject to } (x^\e(\cdot),u(\cdot)) \text{ admissible }.
\earray \eeq To find the minimum of $J(x,\al, u\cd, \lambda)$, it
suffices to select $u\cd$ such that $E(x^\e(T)+\lambda-z)^2$
is minimized. We regard this part as $J^\e(x,\al,u\cd)$ in the
sequel.
Let
$v^\e(x,\al)=\inf_{u(\cdot)}J^\e(x,\al,u(\cdot))$ be the value
function. First define
\beq{r}%
\barray \rho(t,i)=B(t,i)[\sigma(t,i)\sigma'(t,i)]^{-1}B'(t,i), i\in
\{1,2,\ldots,m\}. \earray \eeq
 Consider the following
two systems of ODEs for $i=1,2,\ldots, m$: \beq{p}\barray
\ad\dot{P}^\e(t,i)=P^\e(t,i)[\rho(t,i)-2r(t,i)]-\sum^m_{j=1}q^\e_{ij}P^\e(t,j)\\
\ad P^\e(T,i)=1 .\earray \eeq and
 \beq{h}\barray
\dot{H}^\e(t,i)=\ad H^\e(t,i)r(t,i)-\frac{1}{P^\e(t,i)}\sum^m_{j=1}q^\e_{ij}P^\e(t,j)H^\e(t,j)\\
\ad+\frac{H^\e(t,i)}{P^\e(t,i)}\sum^m_{j=1}q^\e_{ij}P^\e(t,j),\\
H^\e(T,i)=\ad 1. \earray \eeq
The existence and uniqueness of solutions to the above two systems
of equations are easy to obtain since
they are both linear
in the continuous state variable.
Applying the generalized It\^{o}'s formula to
\bea%
v^\e(t,x^\e(t),i)=P^\e(t,i)(x^\e(t)+(\lambda-z)H^\e(t,i))^2, \eea by
employing the completing square techniques, we obtain
\beq{eq1}\barray
\ad\! d P^\e(t,i)[x^\e(t)+(\lambda-z)H^\e(t,i)]^2\\
\aad \
=2P^\e(t,i)[x^\e(t)+(\lambda-z)H^\e(t,i)]dx^\e(t)+P^\e(t,i)(dx^\e(t))^2
\\
\aad \quad +\sum^m_{j=1}q^\e_{ij}P^\e(t,j)
[x^\e(t)+(\lambda-z)H^\e(t,j)]^2dt\\
\aad \quad
+\dot{P}^\e(t,i)[x^\e(t)+(\lambda-z)H^\e(t,i)]^2dt+2P^\e(t,i)[x^\e(t)
+(\lambda-z)H^\e(t,i)](\lambda-z)\dot{H}^\e(t,i)dt.\earray\eeq
Therefore, after plugging in the dynamic equation of $ P^\e(t,i)$ and
$H^\e(t,i)$,  integrating from $0$ to $T$,
and taking expectation, we obtain
\beq{ex}%
\barray%
\ad E[x^\e(T)+\lambda-z]^2\\
\aad \ = P^\e(0,\al)[x+(\lambda-z)H^\e(0,\al)]^2\\
\aad\quad +E\int^T_0(\lambda-z)^2\sum^m_{j=1}q^\e_{ij}P^\e(t,j)[H^\e(t,j)-H^\e(t,i)]^2dt\\
\aad\quad
+E\int^T_0P^\e(t,i)(u(t)-u^{\e,*}(t))'(\sigma(t,i)\sigma'(t,i))(u(t)-u^{\e,*}(t))dt.
\earray%
\eeq%
This leads to the optimal control of the form
\beq{op}
u^{\e,*}(t,\al^\e(t),x^\e(t))=-(\sigma(t,\al^\e(t))\sigma'(t,\al^\e(t)))^{-1}
B'(t,\al^\e(t))[x^\e(t)+(\lambda-z)H^\e(t,\al^\e(t))].\eeq
To proceed, we state a lemma below for subsequent use.
The proof of the lemma is omitted.

\begin{lem}\label{bdd}
The following assertions hold.
\begin{itemize}
\item
  The solutions of equations \eqref{p} and \eqref{h} satisfy
  $0<P^\e(t,i)\le c$ and $0<H^\e(t,i)\le 1$ for all $t\in
  [0,T],i=1,2,\ldots,m$.
  \item
For $i\in\mathcal{M}$, the solutions of \eqref{p} and \eqref{h} are
uniformly Lipschitz on $[0,T]$.
\end{itemize}
\end{lem}

\section{Weak Convergence Results}\label{sec:weak}

Although the optimal solution
of the mean-variance control problem
for the  regime-switching system can be obtained using the methods developed in
\cite{ZhouY}, the difficulty is that $|\M|=m$ is large and we
have to solve a large-scale system, which is computationally
intensive and practically unattractive.
As a viable alternative, we focus on an decomposition-aggregation approach.

Assume that $\wdt Q$ is of the block-diagonal form
$\wdt Q=\diag(\wdt Q^1,\ldots, \wdt Q^l)$ in which $\wdt Q^k\in
\mathbb{R}^{m_k \times m_k}$ are irreducible for $k=1,2,\ldots,l$
and $\sum^l_{k=1}m_k=m $, and $\wdt Q^k$ denotes the $kth$ block matrix
in $\wdt Q$. Let $\mathcal{M}_k=\{s_{k1},s_{k2},\ldots,s_{km_k}\}$
denote the states corresponding to $\wdt Q^k$ and let
$\mathcal{M}=\mathcal{M}_1\cup\mathcal{M}_2\ldots
\cup\mathcal{M}_l=\{s_{11},s_{12},\ldots,s_{1m_1},\ldots,s_{l1},s_{l2},\ldots,s_{lm_l}
\}.$
The slow and fast components are coupled through weak and strong
interactions in the sense that the underlying Markov chain
fluctuates rapidly within a single group $\mathcal{M}_k$ and jumps
less frequently among groups $\mathcal{M}_k$ and $\mathcal{M}_j$ for
$k\neq j$.

By aggregating the states in $\mathcal{M}_k$ as one state $k$, we
can obtain an aggregated process $\lbar{\al}^\e(\cdot)$. That is,
$\lbar{\al}^\e(t)=k$ when $\al^\e(t)\in \mathcal{M}_k$. By virtue of
\cite[Theorem7.4]{Yin&Zhang},  $\lbar{\al}^\epsilon(\cdot)$
converges weakly to $\lbar{\al}(\cdot)$ whose generator is
given by
\beq{barQ-def}
\overline{Q}=\diag(\mu^1,\mu^2,\ldots,\mu^l)
\wdh{Q}\diag(\mathbbold{1}_{m_1},\mathbbold{1}_{m_2},\ldots,
\mathbbold{1}_{m_l}),\eeq where $\mu^k$ is the stationary
distribution of $\wdt Q^k,k=1,2,\ldots,l$, and
$\mathbbold{1}_n=(1,1,\ldots,1)\in \mathbb{R}^n$.
 Define an operator $\mathcal{L}^\e$ by
\beq{L-e}%
\barray%
 \mathcal{L}^\e f(x,t,\iota)=\ad
{\partial f(x,t,\iota)\over \partial t} +
 [r(t,\iota)x
 +B(t,\iota)u(t)]\frac{\partial f(x,t,\iota)}{\partial
 x}\\[.3cm]
 \ad+\frac{1}{2}[u'(t)\sigma(t,\iota)\sigma'(t,\iota)
 u(t)]\frac{\partial^2
 f(x,t,\iota)}{\partial x^2}  +Q^\e f(x,t,\iota), \ \iota \in \M,
\earray
\eeq
where
\beq{mar1}%
 Q^\e f(x,t,\cdot)(\iota)=\sum_{\ell\neq
\iota }q^\e_{ \iota \ell}(f(x,t,\ell)-f(x,t,\iota)),\eeq
and for each $\iota \in \M$, $f(\cdot,\cdot,\iota) \in C^{2,1}$
(that is, $f\cd$ has continuous derivatives up to the second order
with respect to $x$ and continuous derivative with respect to $t$ up
to the first order). Note that the operator, in fact, is $u$-dependent,
so it may be written as $\op^u$.
In this paper, we work with a fixed $u$. We could also consider a
feedback system with $\e$-dependence in the control. In such a setup,
we can use a relaxed control formulation. However, we will not
proceed in this line here.
Define
\beq{barL-e}%
\barray%
 \lbar{\mathcal {L}} f(x,t,k)=\ad
{\partial f(x,t,k)\over \partial t} +
 [\lbar r(t,k)x
 +\lbar B(t,k)u(t)]\frac{\partial f(x,t,k)}{\partial
 x}\\[.3cm]
 \aad+\frac{1}{2}[u'(t)\lbar \sigma(t,k)\lbar\sigma'(t,k)
 u(t)]\frac{\partial^2
 f(x,t,k)}{\partial x^2}  +\lbar Q f(x,t,k), \ k \in \lbar\M,
\earray
\eeq
where  $\lbar Q$ is defined in \eqref{barQ-def} and
\bea \ad \lbar{r}(t,k)=\disp
\sum^{m_k}_{j=1}\mu^k_jr(t,s_{kj}), \\
\ad \lbar{B}(t,k)
=\sum^{m_k}_{j=1}\mu^k_jB(t,s_{kj}),\\
 \ad \lbar{\sigma}^2(t,k)=\disp
\sum^{m_k}_{j=1}\mu^k_j\sigma^2(t,s_{kj}).\eea
The
following theorems are concerned with the weak convergence of a pair
of processes.

\begin{thm}\label{weak-conv}
Suppose that the martingale problem with operator $\lbar {\cal L}$
defined in {\rm \eqref{barL-e}}
has a unique solution for each initial condition. Then
the pair of processes $(x^\e\cd,\lbar\al^\e\cd)$ converges weakly to
$(x\cd,\lbar\al\cd)$, which is the solution of the martingale problem
with operator $\lbar {\mathcal {L}}$.
\end{thm}

\para {Proof.}
The proof is divided into the following steps. First, we prove the
tightness of $x^\e\cd$. Once the tightness is verified, we proceed
to obtain the convergence by using a martingale problem formulation.

Step (i): Tightness.
 We first show that a priori bound holds.

\begin{lem}\label{xb}
Let $x^\e(t)$ denote  flow of system  corresponding to $\al^\e(t)$.
Then $$\sup_{0\le t\le T}E|x^\e(t)|^2=O(1).$$
\end{lem}

 \para{Proof.} Recall that \bea
 dx^\e(t)\ad=[r(t,\al^\e(t))x^\e(t)
 -\rho(t,\al^\e(t))x^\e(t)-\rho(t,\al^\e(t))(\lambda-z)\lbar{H}(t,
\lbar{\al}^\e(t))]dt\\
 \aad \quad +\sum^d_{i=1}\sqrt{\left(\sum^{d_1}_{n=1}
 u^{\e,*}_n(t,x^\e(t),\al^\e(t))\sigma_{ni}(t,\al^\e(t))\right)^2}dw_i(t)\\
 x^\e(0)\ad=x.\eea
 So,\bea
E|x^\e(t)|^2\ad \le K|x|^2+E\left|\int^t_0(r(\nu,\al^\e
(\nu))+\rho(\nu,\al^\e(\nu)))x^\e(\nu))d\nu\right|^2
\\[0.3cm]
\aad \ +KE\int^t_0\left(\sum^{d_1}_{n=1}
 u^{\e,*}_n(\nu,x^\e(\nu),\al^\e(\nu))\sigma_{ni}(\nu,\al^\e(\nu))\right)^2d\nu
 \\[0.4cm]
\ad \le K+KE\int^t_0|x^\e(\nu)|^2d\nu.\eea Here, recall that
$\sigma(t,\al^\e(t))=(\sigma_{ni}(t,\al^\e(t)))\in \rr^{d_1\times
d}$ and note that $u^{\e,*}_n$ is the $n$th component of the $d_1$
dimensional variable. Using properties of stochastic integrals,
H\"{o}lder inequality, and boundedness of
$r(\cdot),B(\cdot),\sigma(\cdot)$,  by Gronwall's inequality, we
obtain the second moment bound of $x^\e(t)$ as desired. \qed

\begin{lem}\label{tgt}
$\{x^\e\cd\}$ is tight in $D([0,T]: \rr)$, the space of
real-valued functions defined on $[0,T]$ that
are right-continuous, and have left limits endowed with the Skorohod
topology.\end{lem}

\para{Proof.}
Denote by  ${\cal F}_t^\e$
the $\sigma$-algebra generated
$\{w(s),\al^\e(s): s\le t\}$ and by
$E^\e_t$ the conditional expectation
w.r.t. ${\cal F}^\e_t$. For any $T<\infty$, any $0\le t\le T$,
any  $s>0$, and any $\delta>0$ with $0< s\le \delta$, by
properties of stochastic integral and boundedness of coefficients,
\bea E^\e_t|x^\e(t+s)-x^\e(t)|^2\ad\le
KE^\e_t\int^{t+s}_t|(r(\nu,\al^\e(\nu))+\rho(\nu,\al^\e(\nu)))
x^\e(\nu)|^2d\nu\\
\aad \quad +KE^\e_t\int^{t+s}_t(\sum^{d_1}_{n=1}
 u^{\e,*}_n(\nu,x^\e(\nu),\al^\e(\nu))\sigma_{ni}(\nu,\al^\e(\nu)))^2d\nu\\
\ad \le Ks+KE^\e_t\int^{t+s}_t|x^\e(\nu)|^2d\nu.\eea
Thus
 we
have
 $$ \lim_{\delta\to 0}\limsup_{\e
\to  0} \sup_{0\le s\le \delta}
 \Big\{ E [ E^\e_t  |x^\e(t+s)-x^\e(t)|^2] \Big\}=0.$$
 Then the tightness criterion
\cite[Theorem 3]{Kush} yields that process $x^\e(\cdot)$ is tight.
Here and in the following part, $K$ is a generic constant which
takes different values in different context.

Step (ii): Using the techniques given in \cite[Lemma
7.18]{Yin&Zhang}, it can be shown that the martingale problem with
operator $\lbar{\mathcal{L}} $ has a unique solution for each
initial condition.

Step (iii): To complete the proof, we characterize the limit process.
Since $(x^\e\cd,\lbar\al^\e\cd)$ is tight, we can
 extract a weakly convergent
subsequence. For notional simplicity, we still denote the
subsequence by $(x^\e\cd,\lbar\al^\e\cd)$ with limit
$(x\cd,\lbar\al\cd)$. By Skorohod representation with no change of
notation, we may assume $(x^\e\cd,\lbar\al^\e\cd)$ converges to
$(x\cd,\lbar\al\cd)$ w.p.1. We next show that the limit
$(x\cd,\lbar\al\cd)$ is a solution of the martingale problem with
operator $\lbar{\mathcal{L}}$ defined by  \eqref{barL-e}.

\begin{lem}\label{mar2} The process $x(\cdot)$ is the solution of the
martingale problem with the operator $\lbar{\mathcal{L}}$.
\end{lem}

\para {Proof.}
To obtain the desirable result, we need to show
  $$f(x(t),t,\lbar{\alpha}(t))-f(x,0,\al)-\int^t_0\lbar{\mathcal{L}}
  f(x(\nu),\nu,
  \lbar{\alpha}(\nu))d\nu
\ \hbox{ is a martingale,}
$$
for $i \in\mathcal{M}, f(\cdot,i)\in C^{2,1}_0([0,T],\rr^r).$ This can
be done by showing that for any integer $n>0$, any bounded and
measurable function $h_p(\cdot,\cdot)$ with $p\le n$, and  any
$t,s,t_p>0$ with $t_p\le t <t+s\le T$,
\bea \ad
E\prod^n_{p=1}h_p(x^\e(t_p),\lbar{\al}^\e(t_p))[f(x(t+s),
t+s,\lbar{\alpha}(t+s))
-f(x(t),t,\lbar{\alpha}(t))\\
\aad \hspace*{0.8in}
-\int^{t+s}_t\lbar{\mathcal{L}}f(x(\nu),\nu,
\lbar{\alpha}(\nu))d\nu]=0.\eea
We further deduce that
\beq{eq:lim-1}\barray\ad  \lim_{\e\to
0}E\prod^n_{p=1}h_p(x^\e(t_p),\lbar{\al}^\e(t_p))(f(x^\e(t+s),t+s,
\lbar{\alpha}^\e
(t+s))-f(x^\e(t),t,\lbar{\alpha}^\e(t))
\\
\aad \
=E\prod^n_{p=1}h_p(x(t_p),\lbar{\al}(t_p))(f(x(t+s),t+s,\lbar{\alpha}(t+s))
-f(x(t),t,\lbar{\alpha}(t)).\earray\eeq
Moreover,
\beq{eq:lim-2}\barray\ad  \lim_{\e\to
0}E\prod^n_{p=1}h_p(x^\e(t_p),\lbar{\al}^\e(t_p))
\Big[ \int^{t+s}_t {\partial f(x^\e(\nu),\nu,\lbar \al^\e(\nu))
\over \partial \nu} d \nu\Big] \\
\aad =
E\prod^n_{p=1}h_p(x(t_p),\lbar{\al}(t_p))
\Big[ \int^{t+s}_t {\partial f(x(\nu),\nu,\lbar \al(\nu))
\over \partial \nu} d \nu\Big]\earray\eeq
by the weak convergence of $(x^\e\cd, \lbar\al^\e\cd)$ and the
Skorohod representation.

For any $f(\cdot)$ chosen above, define
$$\wdh{f}(x^\e(t),t,\al^\e(t))=\sum^l_{i=1}
f(x^\e(t),t,i)I_{\{\alpha^\e(t) \in \mathcal{M}_i\}}$$ since
$(x^\e(t), \al^\e(t))$ is a Markov process, we have
 \bea \wdh{f}(x
^\e(t),t,\al^\e(t))-\wdh{f}(x,0,\alpha)-\int^t_0 \mathcal{L}^\e
\wdh{f}(x^\e(\nu),\nu,\al^\e(\nu))d\nu \eea is a martingale.
 Consequently,
\bea \ad E\disp
\prod^n_{p=1}h_p(x^\e(t_p),
\lbar{\al}^\e(t_p))(\wdh{f}(x^\e(t+s),t+s,\alpha^\e(t+s))
-\wdh{f}(x^\e(t),t,\alpha^\e(t))\\
\aad \quad -\int^{t+s}_t\mathcal{L}^\e
\wdh{f}(x^\e(\nu),\nu,\al^\e(\nu))d\nu)=0.\eea
 Note that
 $\wdh{f}(x^\e(t),t,\al^\e(t))
 =f(x^\e(t),t,\lbar{\alpha}^\e(t))$.

Next we need to show that \bea \ad
\lim_{\e \to
0}E\prod^n_{p=1}h_p(x^\e(t_p),\lbar{\al}^\e(t_p))\int^{t+s}_t\mathcal{L}^\e
\wdh{f}(x^\e(\nu),\nu, \alpha^\e (\nu))d\nu\\
\aad \quad =
E\prod^n_{p=1}h_p(x(t_p),\lbar{\al}(t_p))\int^{t+s}_t\lbar{\mathcal{L}}
f(x(\nu),\nu,\lbar{\alpha}(\nu))d\nu  .\eea
 Note that we can rewrite
 $E\prod^n_{p=1}h_p(x^\e(t_p),\lbar{\al}^\e(t_p))\int^{t+s}_t\mathcal{L}^\e
\wdh{f}(x^\e(\nu),\nu,\al^\e(\nu))d\nu$ as  \bea  \ad
E\prod^n_{p=1}h_p(x^\e(t_p),\lbar{\al}^\e(t_p))
[\int^{t+s}_t\sum^l_{k=1}\sum^{m_k}_{j=1}
Q^\e\wdh{f}(x^\e(\nu),\nu,\cdot)(s_{kj})I_{\{\alpha^\e(\nu)=s_{kj}\}}d\nu\\
\aad\ +\int^{t+s}_t\sum^l_{k=1}\sum^{m_k}_{j=1}
\frac{\partial{\wdh{f}}(x^\e(\nu),\nu,s_{kj})}{\partial
x}I_{\{\alpha^\e(\nu)=s_{kj}\}}[
r(\nu,s_{kj})x^\e(\nu)+B(\nu,s_{kj})u(\nu)]d\nu\\
\aad \ +\int^{t+s}_t\frac{1}{2}
\sum^l_{k=1}\sum^{m_k}_{j=1}[u'(\nu)
\sigma(\nu,s_{kj})\sigma'(\nu,s_{kj})u(\nu)]\frac{\partial^2
 \wdh{f}(x^\e(\nu),
 \nu,s_{kj})}{\partial x^2}I_{\{\alpha^\e(\nu)=s_{kj}\}}]d\nu.\eea
 Since $\wdt Q^k\mathbbold{1}_{m_k}=0$, we
have
$$Q^\e\wdh{f}(x^\e(t),t,\cdot)(s_{kj})=\wdh{Q}\wdh{f}(x^\e(t),t,\cdot)(s_{kj}).$$
We decompose
$$E\prod^n_{p=1}h_p(x^\e(t_p),\lbar{\al}^\e(t_p))\int^{t+s}_t\mathcal{L}^\e
\wdh{f}(x^\e(\nu),\nu,\al^\e(\nu))d\nu$$ as
$H^\e_1(t+s,t)+H^\e_2(t+s,t).$ In which
\bea  H^\e_1(t+s,t)\ad =
E\prod^n_{p=1}h_p(x^\e(t_p),\lbar{\al}^\e(t_p))\\
\aad \ \times\Bigg[\sum^l_{k=1}\sum^{m_k}_{j=1}\int^{t+s}_t \mu^k_j
\frac{\partial{\wdh{f}}(x^\e(\nu),\nu,s_{kj})}{\partial
x}I_{\{\lbar{\alpha}^\e(\nu)=k\}}[r(\nu,s_{kj})
x^\e(\nu) +B(\nu,s_{kj})u(\nu)]d\nu\\
\aad \quad +\frac{1}{2}\sum^l_{k=1}\sum^{m_k}_{j=1}\int^{t+s}_t
\mu^k_j[u'(\nu)\sigma(\nu,s_{kj})\sigma'(\nu,s_{kj})u(\nu)]\frac{\partial^2
 \wdh{f}(x^\e(\nu),\nu,s_{kj})}
 {\partial x^2}I_{\{\lbar{\alpha}^\e(\nu)=k\}}d\nu\\
 \aad\quad +\sum^l_{k=1}\sum^{m_k}_{j=1}\int^{t+s}_t
\mu^k_j\wdh{Q}\wdh{f}(x^\e(\nu),\nu,\cdot)(s_{kj})
I_{\{\lbar{\alpha}^\e(\nu)=k\}}d\nu\Bigg]\eea and $H^\e_2(t+s,t)$
can be represented as \bea
\ad H^\e_2(t+s,t)\\
 \ad =E\prod^n_{p=1}h_p(x^\e(t_p),\lbar{\al}^\e(t_p)
 \Big(\sum^l_{k=1}\sum^{m_k}_{j=1}\int^{t+s}_t
(I_{\{\alpha^\e(\nu)=s_{kj}\}}-\mu^k_jI_{\{\lbar{\alpha}^\e(\nu)=k\}})
\frac{\partial{\wdh{f}}(x^\e(\nu),\nu,s_{kj})}{\partial
x}\times\\
\aad\ [r(\nu,s_{kj})x^\e(\nu)+B(\nu,s_{kj})u(\nu)]d\nu+\sum^l_{k=1}
\sum^{m_k}_{j=1}\int^{t+s}_t (I_{\{\alpha^\e(\nu)=s_{kj}\}}-\mu^k_j
I_{\{\lbar{\alpha}^\e(\nu)=k\}})\wdh{Q}\times\\
\aad\
\wdh{f}(x^\e(\nu),\nu,\cdot)(s_{kj})d\nu+\frac{1}{2}\sum^l_{k=1}
\sum^{m_k}_{j=1}\int^{t+s}_t
(I_{\{\alpha^\e(\nu)=s_{kj}\}}-\mu^k_jI_{\{\lbar{\alpha}^\e(\nu)=k\}})\times\\
\aad\
[u'(\nu)\sigma(\nu,s_{kj})\sigma'(\nu,s_{kj})u(\nu)]\frac{\partial^2
 \wdh{f}(x^\e(\nu),\nu,s_{kj})}
 {\partial x^2}d\nu\Big).\eea
 By virtue of \lemref{re}, \cite[Theorem7.14]{Yin&Zhang},
Cauchy-Schwartz inequality, boundedness of $h_p(\cdot)$, $r(\cdot)$
and $B(\cdot)$, for each $k=1,2,\ldots, l; j=1,2,\ldots, m_k$, as
$\e \to 0 $ \bea \ad
E|\prod^n_{p=1}h_p(x^\e(t_p),\lbar{\al}^\e(t_p))\int^{t+s}_t
(I_{\{\alpha^\e(\nu)=s_{kj}\}}-\mu^k_jI_{\{\lbar{\alpha}^\e(\nu)=k\}})
\frac{\partial{\wdh{f}}(x^\e(\nu),\nu,s_{kj})}{\partial x} \\
\aad\qquad \times
[r(\nu,s_{kj})x^\e(\nu)+B(\nu,s_{kj})u(\nu)]d\nu|^2\to 0. \eea
Similarly
 as $\e \to 0$, \bea \ad
E|\prod^n_{p=1}h_p(x^\e(t_p),\lbar{\al}^\e(t_p))\int^{t+s}_t
(I_{\{\alpha^\e(\nu)=s_{kj}\}}-\mu^k_jI_{\{\lbar{\alpha}^\e(\nu)=k\}})\\
\aad \quad \times
[u'(\nu)\sigma(\nu,s_{kj})\sigma'(\nu,s_{kj})u(\nu)]\frac{\partial^2
 \wdh{f}(x^\e(\nu),\nu,s_{kj})}{\partial x^2}d\nu|^2 \to 0,\eea
 and
$$
E|\prod^n_{p=1}h_p(x^\e(t_p),\lbar{\al}^\e(t_p))\int^{t+s}_t
(I_{\{\alpha^\e(\nu)=s_{kj}\}}-\mu^k_j
I_{\{\lbar{\alpha}^\e(\nu)=k\}})\wdh{Q}\wdh{f}(x^\e(\nu),\nu,
\cdot)(s_{kj})d\nu|^2\to
0.$$ Therefore,
$H^\e_2(t+s,t)$ converges to $0$ in probability.
On the other hand, we obtain
  \beq{h2}\barray \ad
E\prod^n_{p=1}h_p(x^\e(t_p),\lbar{\al}^\e(t_p))\sum^l_{k=1}\sum^{m_k}_{j=1}
\int^{t+s}_t\mu^k_j\frac{\partial{\wdh{f}}(x^\e(\nu),\nu,s_{kj})}
{\partial x}[r(\nu,s_{kj})x^\e(\nu)+B(\nu,s_{kj})u(\nu)]\\
\aad \quad \hfill \times I_{\{\lbar{\alpha}^\e(\nu)=k\}}d\nu\\
\aad \ \to
\sum^l_{k=1}\sum^{m_k}_{j=1}E\prod^n_{p=1}h_p(x(t_p),\lbar{\al}(t_p))
\int^{t+s}_t\mu^k_j\frac{\partial{f}(x(\nu),
\nu,\lbar{\alpha}(\nu))}{\partial
x}[r(\nu,s_{kj})x(\nu)+B(\nu,s_{kj})u(\nu)]\\
\aad \quad \hfill \times I_{\{\lbar{\alpha}(\nu)=k\}}d\nu\\
\aad \ =\sum^l_{k=1}E\prod^n_{p=1}h_p(x(t_p),\lbar{\al}(t_p))
\int^{t+s}_t\frac{\partial{f}(x(\nu),\nu,\lbar{\alpha}(\nu))}{\partial
x}[\lbar{r}(\nu,\lbar{\alpha}(\nu))x(\nu)+\lbar{B}
(\nu,\lbar{\alpha}(\nu))u(\nu))]\\
\aad \quad \hfill \times I_{\{\lbar{\alpha}(\nu)=k\}}d\nu\\
\aad \ =E\prod^n_{p=1}h_p(x(t_p),\lbar{\al}(t_p))\int^{t+s}_t
\frac{\partial{f}(x(\nu),\nu,\lbar{\alpha}(\nu))}{\partial{x}}
[\lbar{r}(\nu,\lbar{\alpha}(\nu))x(\nu)+
\lbar{B}(\nu,\lbar{\alpha}(\nu))u(\nu)]d\nu.\earray\eeq Similarly,
\beq{h3}\barray \ad E\prod^n_{p=1}h_p(x^\e(t_p),\lbar{\al}^\e(t_p))
\sum^l_{k=1}\sum^{m_k}_{j=1}\int^{t+s}_t
\mu^k_j[u'(\nu)\sigma^2(\nu,s_{kj})u(\nu)]\frac{\partial^2
 \wdh{f}(x^\e(\nu),\nu,s_{kj})}
 {\partial x^2}I_{\{\lbar{\alpha}^\e(\nu)=k\}}d\nu\\
\aad \ \to\
E\prod^n_{p=1}h_p(x(t_p),\lbar{\al}(t_p))\int^{t+s}_t\frac{\partial^2
 f(x(\nu),\nu,\lbar{\alpha}(\nu))}{\partial x^2}[u'(\nu)
 \lbar{\sigma}^2(\nu,\lbar{\alpha}(\nu))u(\nu)]d\nu.\earray\eeq
Note that $$ \sum^l_{k=1}\sum^{m_k}_{j=1}\int^{t+s}_t
\mu^k_jI_{\{\lbar{\alpha}^\e(\nu)=k\}}
\wdh{Q}\wdh{f}(x^\e(\nu),\nu,\cdot)(s_{kj})d\nu
=\int^{t+s}_t\lbar{Q}f(x^\e(\nu),\nu,\cdot)(\lbar{\alpha}^\e(\nu))d\nu.$$
So as $\e \to 0$,%
\beq{h22}%
\barray\disp%
\ad
E\prod^n_{p=1}h_p(x^\e(t_p),\lbar{\al}^\e(t_p))\int^{t+s}_t\lbar{Q}f(x^\e(\nu),
\nu,
\cdot)(\lbar{\alpha}^\e(\nu))d\nu \\
\aad \ \to
E\prod^n_{p=1}h_p(x(t_p),\lbar{\al}(t_p))\int^{t+s}_t\lbar{Q}
f(x(\nu),\nu,\cdot)(\lbar{\alpha}(\nu))d\nu.
\earray%
\eeq
Combining the results from \eqref{h2} to \eqref{h22}, we have
 \beq{h1}%
 \barray%
 \disp\ad\lim_{\e \to 0}E\prod^n_{p=1}h_p(x^\e(t_p),\lbar{\al}^\e(t_p))
 \int^{t+s}_t\mathcal{L}^\e
\wdh{f}(x^\e(\nu),\nu,\alpha^\e(\nu))d\nu\\
\aad \ =E\prod^n_{p=1}h_p(x(t_p),\lbar
{\alpha}(t_p))\int^{t+s}_t\lbar{\mathcal{L}}f(x(\nu),\nu,
\lbar{\alpha}(\nu))d\nu\earray \eeq Finally, piecing together the
results obtained, the proof of the theorem is completed.  \qed

To proceed,
we can further deduce the following result. The proof is
omitted.

\begin{thm}\label{4.1}
For $k=1,2,\ldots,l$ and $j=1,2,\ldots,m_k$, $P^\e(t,s_{kj})\to
\overline{P}(t,k)$ and $H^\e(t,s_{kj})\to \overline{H}(t,k)$
uniformly on $[0,T]$ as $\e \to 0$, where $\overline{P}(t,k)$ and
$\overline{H}(t,k)$ are the unique solutions of the following
differential equations for $k=1,2,\ldots,l$, \beq{P2} \barray
\dot{\overline{P}}(t,k)=\ad(\overline{\rho}(t,k)-2\lbar{r}(t,k))\overline{P}(t,k)-
\lbar{Q} \overline{P}(t,\cdot)(k)\\
\overline{P}(T,k)=\ad 1. \earray \eeq and \beq{H2}\barray
\dot{\overline{H}}(t,k)=\ad\overline{r}(t,k)\overline{H}(t,k)-\frac{1}{\overline{P}(t,k)}
\lbar{Q}\overline{P}(t,\cdot)\overline{H}(t,\cdot)(k)\\
\ad+\frac{\overline{H}(t,k)}{\overline{P}(t,k)}\lbar{Q} \overline{P}(t,\cdot)(k)\\
\overline{H}(T,k)=\ad 1. \earray\eeq
\end{thm}

\section{Near Optimality and Numerical Examples}\label{sec:near}
This section establishes  near optimality of the control obtained from the limit system and presents
related numerical results.

\subsection{Near Optimality}
By the convergency of $P^\e(t,s_{kj})$ to $\lbar{P}(t,k)$ and
$H^\e(t,s_{kj})$ to $\lbar{H}(t,k)$, we have $v^\e(t,s_{kj},x)\to
\lbar{v}(t,k,x)$ as $\e \to 0$, in which
$\lbar{v}(t,k,x)=\lbar{P}(t,k)(x+(\lambda-z)\lbar{H}(t,k))^2$. Here,
$\lbar{v}(t,k,x)$ corresponds to the value function of a limit
problem.
In view of \thmref{weak-conv},
for the limit problem,
let $\mathcal{U}$ be the control set
$\mathcal{U}=\{U=(U^1,U^2,\ldots,U^l):
U^k=(u^{k1},u^{k2},\ldots,u^{km_k}),u^{kj}\in
\mathbb{R}^{d_1}\}$. Define \bea \ad
\Phi(t,x,k,U)=\sum^{m_k}_{j=1}\mu^k_jr(t,s_{kj})x
+\sum^{m_k}_{j=1}\mu^k_jB(t,s_{kj})u^{kj}(t) \ \hbox{
and }\\
\ad  \Psi(t,k,U)=((g_1(t,k,U)),\ldots, g_d(t,k,U)) \text{ with }\\
\ad \Psi_i(t,k,U)=\sqrt{\sum^{m_k}_{j=1}{\mu^k_j}
\left(\sum^{d_1}_{n=1}u^{kj}_n\sigma_{ni}(t,s_{kj})\right)^2}.\eea
Here, recall that $\sigma(t,\al^\e(t))=(\sigma_{ni}(t,s_{kj}))\in
\rr^{d_1\times d}$ and note that $u^{kj}_n$ is the $n$th component
of the $d_1$ dimensional variable. The corresponding dynamic system
is given by
 \beq{li} dx(t)=\Phi(t,x(t),\lbar{\al}(t),
U(t))dt+\sum^d_{i=1}\Psi_i(t,\lbar{\al}(t),U(t))dw_i(t),\eeq
where $\overline{\al}(\cdot)\in \{1,2,\ldots,l\}$ is a Markov chain
generated by $\lbar{Q}$ with $\lbar{\al}(0)=\al$.
It can be shown
that
 the optimal control
for this limit problem is \bea
U^*(t,x)=(U^{1*}(t,x),U^{2*}(t,x),\ldots,U^{l*}(t,x))\eea with
\[U^{k*}(t,x)=(u^{k1*}(t,x),u^{k2*}(t,x),\ldots,u^{km_k*}(t,x))\] and
\[u^{kj*}(t,x)=-(\sigma(t,s_{kj})\sigma'(t,s_{kj}))^{-1}B'(t,s_{kj})[x+(\lambda-z)
\lbar{H}(t,k)].\] Using such controls, we construct \beq{con}
u^\e(t,\al^\e(t),
x)=\sum^l_{k=1}\sum^{m_k}_{j=1}I_{\{\al^\e(t)=s_{kj}\}}
u^{kj*}(t,x)\eeq for the original problem. This control can also be
written as if $\al^\e(t)\in\mathcal{M}_k,u^\e(t,\al^\e(t),
x)=-(\sigma(t,\al^\e(t))\sigma'(t,\al^\e(t)))^{-1}
B'(t,\al^\e(t))[x+(\lambda-z)
\lbar{H}(t,\lbar{\al}^\e(t))]$.
It can be shown that our constructed control
is nearly optimal. We present the
following lemmas first.

\begin{lem}\label{re}
For any $k=1,2,\ldots,l,j=1,2,\ldots, m_k$, we have the following
result hold. \beq{ind1}
E\l\int^t_0[I_{\{\al^\e(\nu)=s_{kj}\}}-\mu^k_j
I_{\{\lbar{\al}^\e(\nu)=k\}}]x^\e(\nu)r(\nu,\al^\e(\nu))d\nu\r^2\to
0 \text{ as } \e \to 0.\eeq
\end{lem}

\para{Proof.} For $0<\delta<1$ and any $t\in[0,T]$,
let $N=[t/\e^{1-\delta}]$ partition $[0,t]$ into subintervals of
equal length $\e^{1-\delta}$ and denote the partition boundaries by
$t_k=k\e^{1-\delta}$ for $0\le k\le N-1 $. Define the auxiliary
function \bea \wdt F(\nu)=r(\nu,\alpha^\e(\nu))x^\e(t_k), \text{
for }u\in[t_k,t_{k+1}].\eea \lemref{tgt} shows $$
E|x^\e(t)-x^\e(t_k)|^2=O(\e^{1-\delta})\to 0\text{ as }\e\to 0.$$
for $t\in[t_k,t_{k+1}], 0\le k\le N-1$. Then%
\beq{es}\barray \ad
E\left|\int^t_0[I_{\{\alpha^\e(\nu)=s_{kj}\}}-\mu^k_j
I_{\{\lbar{\alpha}^\e(\nu)=k\}}]x^\e(\nu)r(\nu,\alpha^\e(\nu))d\nu\right|^2\\
\aad\ \le 2E\left|\int^t_0[I_{\{\alpha^\e(\nu)=s_{kj}\}}-\mu^k_j
I_{\{\lbar{\alpha}^\e(\nu)=k\}}]\wdt F(\nu)d\nu\right|^2\\
\aad\qquad +2E\left|\int^t_0[I_{\{\alpha^\e(\nu)=s_{kj}\}}-\mu^k_j
I_{\{\lbar{\alpha}^\e(\nu)=k\}}](x^\e(\nu)r(\nu,\alpha^\e(\nu))-\wdt F(\nu))d\nu\right|^2.\earray\eeq
First, we estimate the last term of \eqref{es}. According to
Cauchy-Schwartz inequality, we have\bea \ad
E\left|\int^t_0[I_{\{\alpha^\e(\nu)=s_{kj}\}}-\mu^k_j
I_{\{\lbar{\alpha}^\e(\nu)=k\}}](x^\e(\nu)r(\nu,\alpha^\e(\nu))-\wdt F(\nu))d\nu\right|^2\\
 \aad\ \le K\int^t_0E|x^\e(\nu)r(\nu,\alpha^\e(\nu))-\wdt F(\nu)|^2d\nu\\
\aad\ \le
K\sum^{N-1}_{k=0}\int^{t_{k+1}}_{t_k}E|x^\e(\nu)r(\nu,\alpha^\e(\nu))-r(\nu,\alpha^\e(\nu))x^\e(t_k)|^2d\nu\\
\aad\ \le KO(\e^{1-\delta})\to 0 \text{ as }\e\to 0.\eea For the
first term of \eqref{es}, for each $k=1,2,\ldots,l$ and
$j=1,2,\ldots,m_k$ define $$
\eta^\e(t)=E\left|\int^t_0[I_{\{\alpha^\e(\nu)=s_{kj}\}}-\mu^k_j
I_{\{\lbar{\alpha}^\e(\nu)=k\}}]\wdt F(\nu)d\nu\right|^2.$$ With the
similar idea involved in \cite[Lemma 7.14]{Yin&Zhang}, we get $$
\disp \sup_{0\le t\le T}\eta^\e(t)=\sup_{0\le t\le
T}\int^t_0O(\e^{1-\delta})d\nu=O(\e^{1-\delta})\to 0 \text{ as }\e
\to 0.$$ Thus, we conclude the proof by combining the above two
parts. \qed

\begin{lem}\label{le2}
For any $k=1,2,\ldots,l,j=1,2,\ldots, m_k$, we have the following
result hold. \beq{ind1a}
E(I_{\{\lbar{\al}^\e(s)=k\}}-I_{\{\lbar{\al}(s)=k\}})^2\to 0 \text{
as } \e \to 0.\eeq
\end{lem}

\para{Proof.} Similar to \cite[Theorem 7.30]{Yin&Zhang}, we can see that
$(I_{\{\lbar{\al}^\e\cd=1\}},\ldots, I_{\{\lbar{\al}^\e\cd=l\}})$
converges weakly to
$(I_{\{\lbar{\al}\cd=1\}},\ldots,I_{\{\lbar{\al}\cd=l\}})$ in $D[0,T]$
as $\e \to 0$. By means of Cram\'{e}r-Word's device, for each $i\in
\mathcal{M}$, $I_{\{\lbar{\al}^\e\cd=i\}}$ converges weakly to
$I_{\{\lbar{\al}\cd=i\}}$. Then the Skorohod representation (with a
little bit of abuse of notation), we may assume
$I_{\{\lbar{\al}^\e\cd=i\}}\to I_{\{\lbar{\al}\cd=i\}}$ w.p.1 without
change of notation. Now by dominance convergence theorem, we can
conclude the proof.\qed

The following result was obtained in \cite{YYWZ}. We state the result
and omit the proof.

 \begin{thm}\label{1}
The control $u^\e(t)$ defined in \eqref{con} is nearly optimal in
that
 $$\lim_{\e \to 0} |J^\e(\al,x,u^\e(\cdot))-v^\e(\al,x)|=0.$$
\end{thm}

Next, we consider the case in which the Markov chain has
transient states. We assume \bea \wdt Q= \left( \begin{array}{l}
 \wdt Q_r \;\,0 \\
 \wdt Q_0 \;\wdt Q_*  \\
 \end{array} \right)
 \eea
 where
 $\wdt Q_r=\diag(\wdt Q^1,\wdt Q^2,\ldots,\wdt Q^l)$,
 $\wdt Q_0=(\wdt Q^1_*,\ldots,\wdt Q^l_*).$ For
 each $k=1,2,\ldots,l$, $\wdt Q^k$ is a generator with dimension
 $m_k\times m_k$, $\wdt Q_*\in \mathbb{R}^{m_*\times m_*}$, $\wdt Q^k_* \in \mathbb{R}^{m_*\times m_k}$, and
 $m_1+m_2+\cdots+
 m_*=m$. The state space of the underlying Markov chain is given by
 $\mathcal{M}=\mathcal{M}_1\cup\mathcal{M}_2\cup
 \ldots\cup\mathcal{M}_*=\{s_{11},\ldots,s_{1m_1},\ldots,s_{l1}\ldots,s_{lm_l},s_{*1},\ldots,s_{*m_*}
 \}$, where $\mathcal{M}_*=\{s_{*1},s_{*2},\ldots,s_{*m_*}\}$
 consists of the transient states. Suppose for $k=1,2,\ldots,l$,
 $\wdt Q^k$ are irreducible, and $\wdt Q_*$ is Hurwitz, i.e., it has
 eigenvalues with negative real parts.
 Let \bea
 \wdh{Q}=
\left( \begin{array}{l}
 \wdh Q^{11} \;\wdh Q^{12}  \\
 \wdh Q^{21} \;\wdh Q^{22}  \\
 \end{array} \right)
\eea where $\wdh{Q}^{11}\in \mathbb{R}^{(m-m_*)\times(m-m*)}$,
$\wdh{Q}^{12}\in \mathbb{R}^{(m-m_*)\times m_*}$, $\wdh{Q}^{21}\in
\mathbb{R}^{m_*\times (m-m_*)}$, and $\wdh{Q}^{22}\in
\mathbb{R}^{m_*\times m_*}$. We define%
\bea
\lbar{Q}_*=\diag(\mu^1,\ldots,\mu^l)(\wdh{Q}^{11}\tilde{\mathbbold{1}}+\wdh{Q}^{12}(a_{m_1},a_{m_2},\ldots,a_{m_l}))
\eea with $\tilde{\mathbbold{1}}=\diag(\mathbbold{1}_{m_1},\ldots,
\mathbbold{1}_{m_l})$, $\mathbbold{1}_{m_j}=(1,\ldots, 1)'\in
\mathbb{R}^{m_j}$ and, for $k=1,\ldots,l$,
\bea%
a_{m_k}=(a_{{m_k},1},\ldots,a_{{m_k},m_*})'=-\wdt Q^{-1}_*\wdt Q^k_*\mathbbold{1}_{m_k}.\eea
Let $\xi$ be a random variable uniformly distributed on $[0,1]$ that
is independent of $\al^\e\cd$. For each $j=1,2,\ldots,m_*$, define
an integer-valued random variable $\xi_j$ by \bea%
\xi_j=I_{\{0 \le\xi\le a_{{m_1},j}\}}+2I_{\{a_{{m_1},j} <\xi \le
a_{{m_1},j}+ a_{{m_2},j} \}}+\cdots+lI_{\{ a_{{m_1},j}+\cdots+
a_{{m_{l-1}},j}<\xi\le 1\}}.\eea Now define the aggregated process
$\lbar{\al}^\e\cd$ by \bea%
 \lbar{\al}^\e(t)= \left\{ \begin{array}{l}
 k,\;\text{ if }\al ^\e(t)\in {\cal M}_k, \\
 \xi _j ,\text{ if }\al ^\e (t) = s_{*j.}  \\
 \end{array} \right.
\eea Note the state space of $\lbar{\al}^\e(t)$ is $\overline{\cal
M}=\{1,2,\ldots,l\}$ and $\lbar{\al}^\e\cd\in D[0,T]$. In addition,
$$P(\lbar{\al}^\e(t)=i|\al^\e(t)=s_{*j})=a_{{m_i},j}.$$ Then according to
\cite[Theorem 4.2]{YZB}, $\lbar{\al}^\e(\cdot)$ converges weakly to
$\lbar{\al}(\cdot)$ where $\lbar{\al}(\cdot)\in
\{1,2,\ldots,l\}$ is the Markov chain generated by $\lbar{Q}_*$. The
following two theorems were also proved in \cite{YYWZ}.

\begin{thm}
As $\e \to 0$, we have $P^\e(s,s_{kj})\to \overline{P}(s,k)$ and
$H^\e(s,s_{kj})\to \overline{H}(s,k)$, for $k=1,2,\ldots, l$,
$j=1,2,\ldots,m_k$, $P^\e(s,s_{*j})\to \overline{P}_*(s,j)$ and
$H^\e(s,s_{*j})\to \overline{H}_*(s,j)$, for $j=1,2,\ldots, m_*$
uniformly on $[0,T]$ where \bea
\overline{P}_*(s,j)=a_{{m_1},j}\overline{P}(s,1)+\cdots+a_{{m_l},j}\overline{P}(s,l),\eea
\bea
\overline{H}_*(s,j)=a_{{m_1},j}\overline{H}(s,1)+\cdots+a_{{m_l},j}\overline{H}(s,l)\eea
and $\overline{P}(s,k)$ and $\overline{H}(s,k)$ are the unique
solutions to the following equations. For $k=1,2,\ldots,l$,
\beq{P2T} \barray
\dot{\overline{P}}(t,k)=\ad(\overline{\rho}(t,k)-2\overline{r}(t,k))\overline{P}(t,k)-
\overline{Q}_*\overline{P}(t,\cdot)(k),\\
\overline{P}(T,k)=\ad 1. \earray \eeq And \beq{H2a}\barray
\dot{\overline{H}}(t,k)=\ad\overline{r}(t,k)\overline{H}(t,k)-\frac{1}{\overline{P}(t,k)}
\overline{Q}_*\overline{P}(t,\cdot)\overline{H}(t,\cdot)(k)\\
\ad+\frac{\overline{H}(t,k)}{\overline{P}(t,k)}\overline{Q}_*\overline{P}(t,\cdot)(k)\\
\overline{H}(T,k)=\ad 1. \earray\eeq
\end{thm}

\begin{thm}\label{thm:2}
Construct \beq{con2} u^\e(t,\al^\e(t),
x)=\sum^l_{k=1}\sum^{m_k}_{j=1}I_{\{\al^\e(t)=s_{kj}\}}
u^{kj*}(t,x)+\sum^{m_*}_{j=1}I_{\{\al^\e(t)=s_{*j}\}}u^{*j*}(t,x)\eeq
for the original problem where
$$u^{*j*}(t,x)=-(\sigma(t,s_{*j})\sigma'(t,s_{*j}))^{-1}
B'(t,s_{*j})[x+(\lambda-z) \lbar{H}_*(t,j)].$$ Then control $u^\e(t,
\al^\e(t),x)$ defined in \eqref{con2} is nearly optimal. That is,
 $$\lim_{\e \to 0} |J^\e(\al,x,u^\e\cd)-v^\e(\al,x)|=0.$$
\end{thm}

\subsection{Numerical Examples}
In this section, we present a couple of examples to demonstrate the
performance of our approximation schemes.
First, let us consider that the Markov chain has only
recurrent states.

\begin{exm}\label{6.1} {\rm We consider the networked system in which the Markov chain
$\al^\e(t)\in\mathcal{M}=\{1,2,3,4\}$, $ t\ge 0$ generated by
$Q^\e= \wdt Q/\e + \wdh Q$ with
\bea  \ad \wdt Q = \left( \begin{array}{rrrr}
  -1 & 1 & 0 & 0 \\
 2 & -2 & 0 & 0 \\
 0 & 0 & -1 & 1 \\
 0 & 0 & 3 &-3 \\
 \end{array} \right), \\
 \ad \wdh Q = \left( \begin{array}{rrrr}
  -2 & 0& 1 & 1 \\
 1 & -2 & 1 & 0 \\
 0 &1 &-1 & 0 \\
 1& 2 & 0 & -3 \\
 \end{array} \right)
\eea the condition of coefficients for the flows of our system are
as follows: Consider the following dynamic system model for
$t\in[0,5]$,
$x^\e(0)=0, r(t,1)=.5, r(t,2)=-.1, r(t,3)=.5, r(t,4)=-.1, B(t,1)=1,
B(t,2)=2, B(t,3)=-1, B(t,4)=-2,
\sigma(t,1)=\sigma(t,2)=\sigma(t,3)=\sigma(t,4)=1$. In this case, we
can classify the Markov chain into two recurrent groups $\M_1$ and
$\M_2$. So the corresponding stationary distribution for
$\mathcal{M}_1$  is $\mu^1=\{\frac{2}{3},\frac{1}{3}\}$ and that of
$\mathcal{M}_2$ is $\mu^2=\{\frac{3}{4},\frac{1}{4}\}$. We
discretize the equations with step size $h=.01$. So, in the
corresponding discrete time setting the time horizon is
$T_h=\frac{5}{h}$. Let $x^\e(t)$ be the optimal trajectory and
$\lbar{x}(t)$ be the nearly optimal trajectory with control taking
as $u^\e\cd$. Sample paths of $\al^\e(t)$, trajectories of
$|x^\e(t)-\lbar{x}(t)|$ are given in Figure 1 for $\e=0.1$ and in
Figure 2 for $\e=0.01$. The results below are based
on computations using $100$ sample paths.  Define \bea%
|P^\e-\lbar{P}|=\ad
\frac{1}{T_h}\sum^{T_h}_{j=1}(|P^\e(jh,1)-\lbar{P}(jh,1)|
+|P^\e(jh,2)-\lbar{P}(jh,1)|\\
\ad+|P^\e(jh,3)-\lbar{P}(jh,2)|+|P^\e(jh,4)-\lbar{P}(jh,2)|),\eea
and \bea
\disp|x^\e-\lbar{x}|=\frac{1}{T_h}\sum^{T_h}_{j=1}|x^\e(jh)-\lbar{x}(jh)|.\eea
Then we have the error bounds given in Table $1$ for different
values of $\e$. }\end{exm}
\begin{table}[H]
\centering \caption{Error Bounds for \exmref{6.1}}
\begin{tabular}{cccc}
\hline
$\e$ & $|P^\e-\lbar{P}|$ & $|x^\e-\lbar{x}|$ & $|J^\e-v^\e|$\\
\hline
0.1 & 1.48 & 0.02 & 0.002\\
0.01 & 0.52 & 0.016 & 0.001\\
0.001 & 0.11 & 0.004 & 0.0003\\
\hline
\end{tabular}
\end{table}

\begin{figure}
\includegraphics[width=80mm]{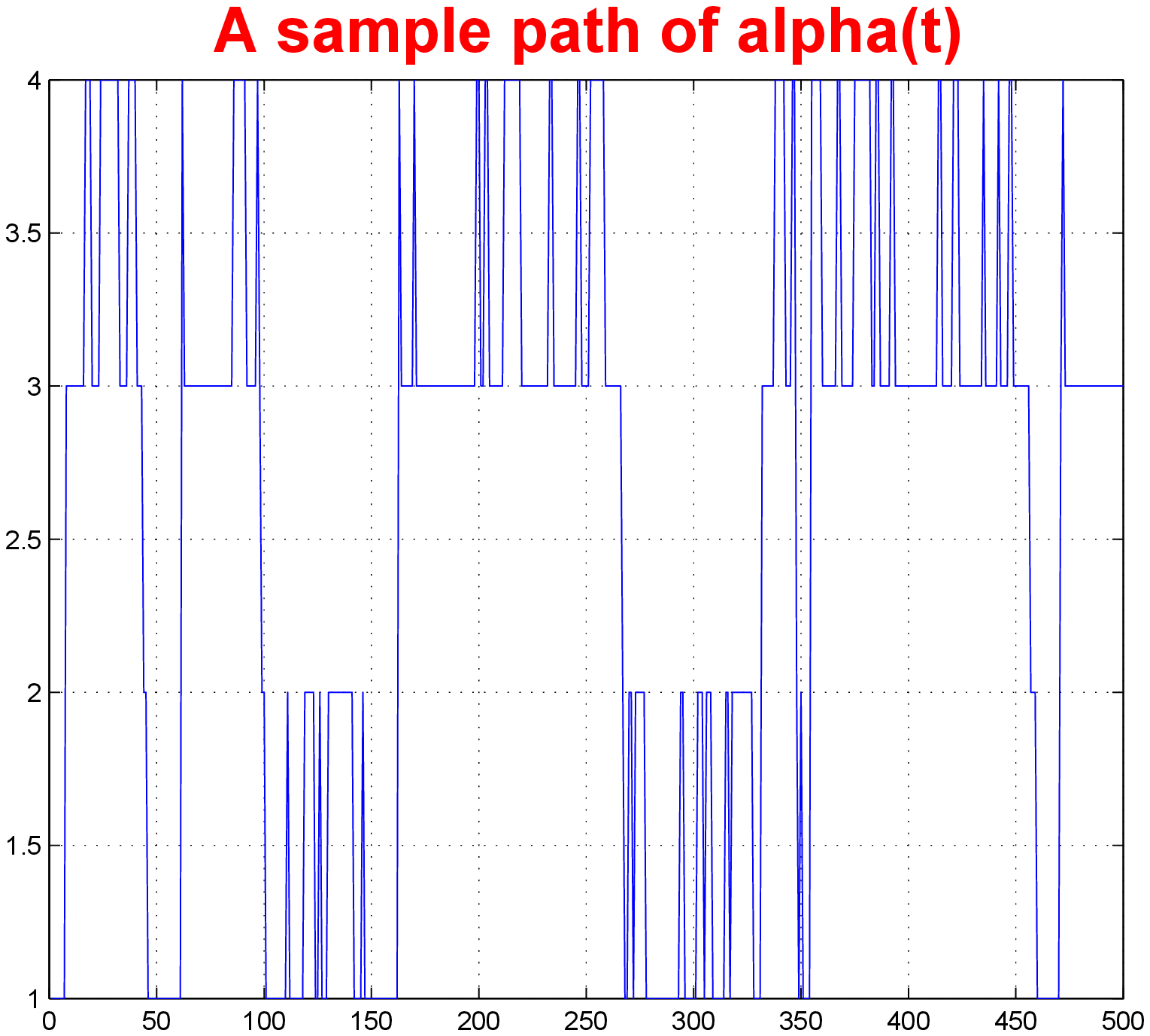}
\includegraphics[width=80mm]{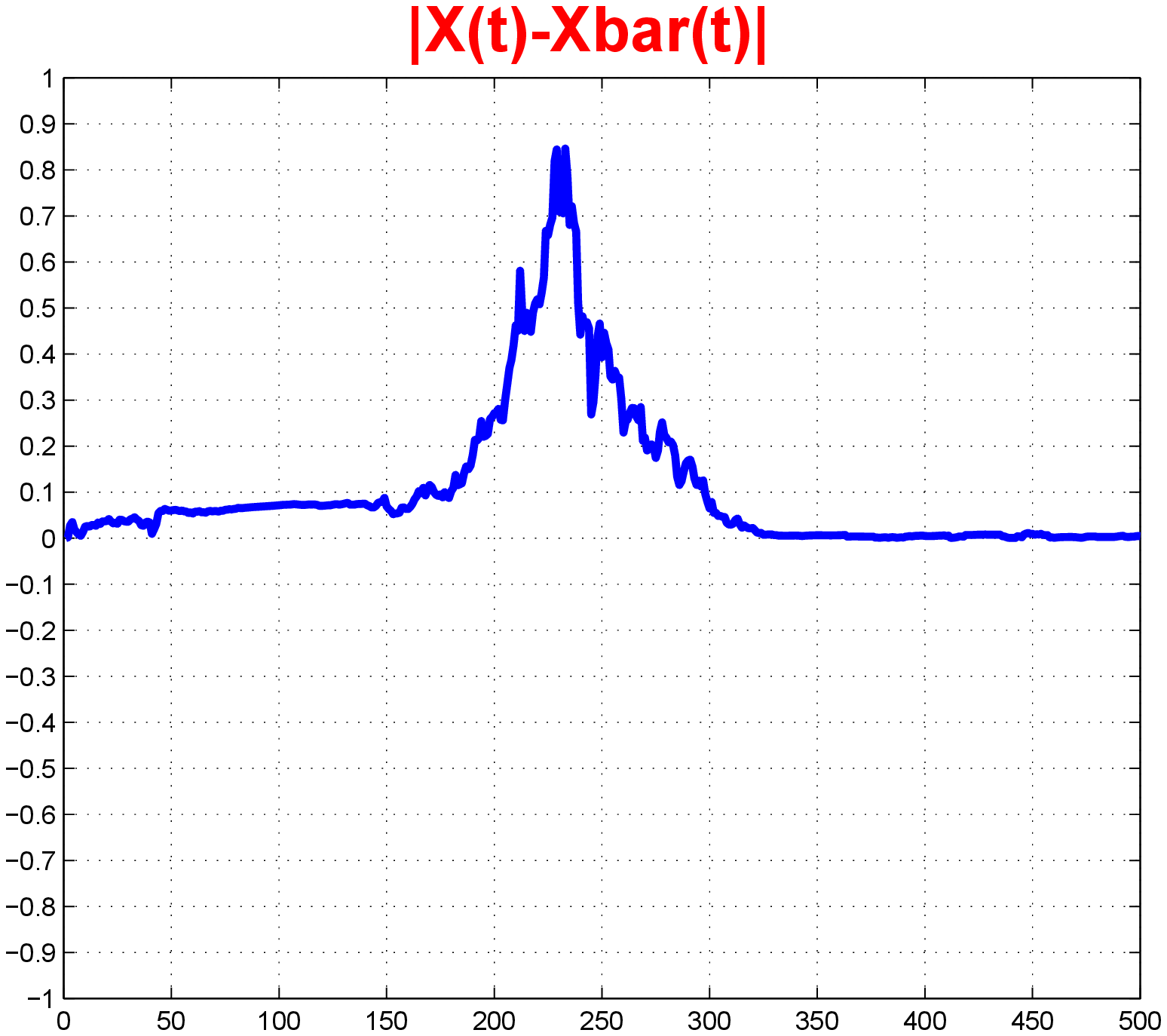}
\caption{Sample path with $\e=0.1$ in \exmref{6.1}}
\end{figure}
\begin{figure}
\includegraphics[width=80mm]{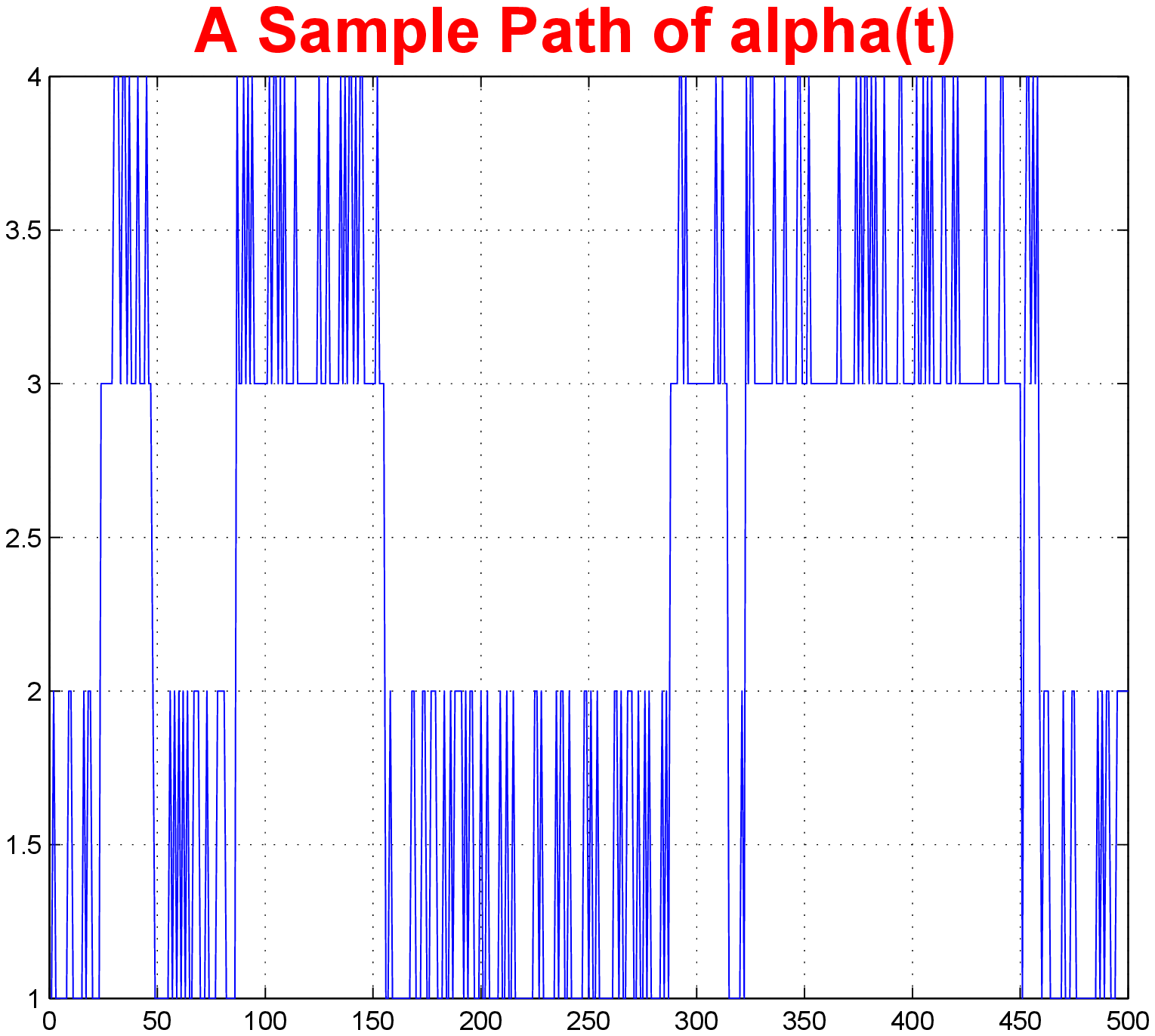}
\includegraphics[width=80mm]{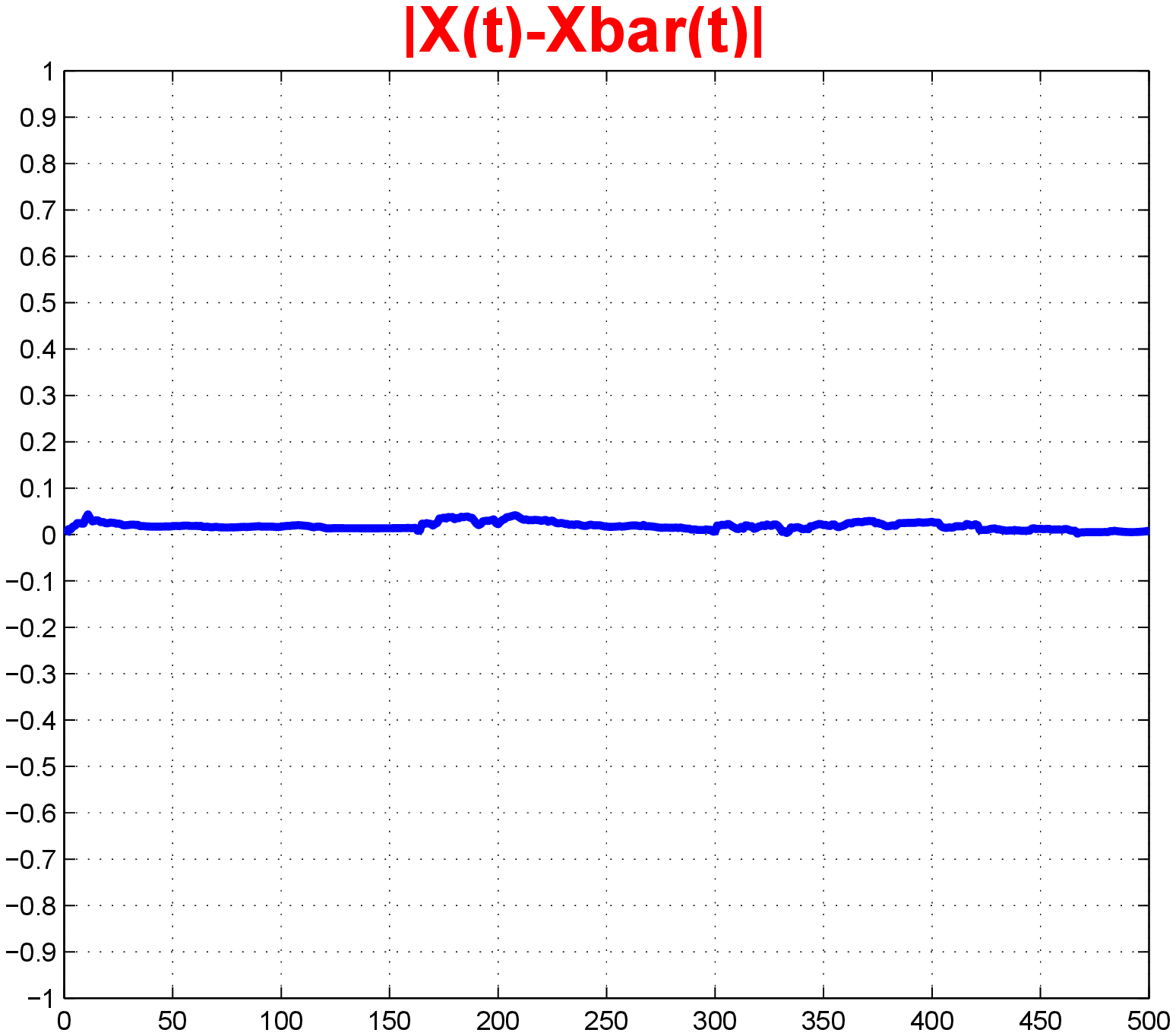}
\caption{Sample path with $\e=0.01$ in \exmref{6.1}}
\end{figure}

\begin{figure}
\includegraphics[width=80mm]{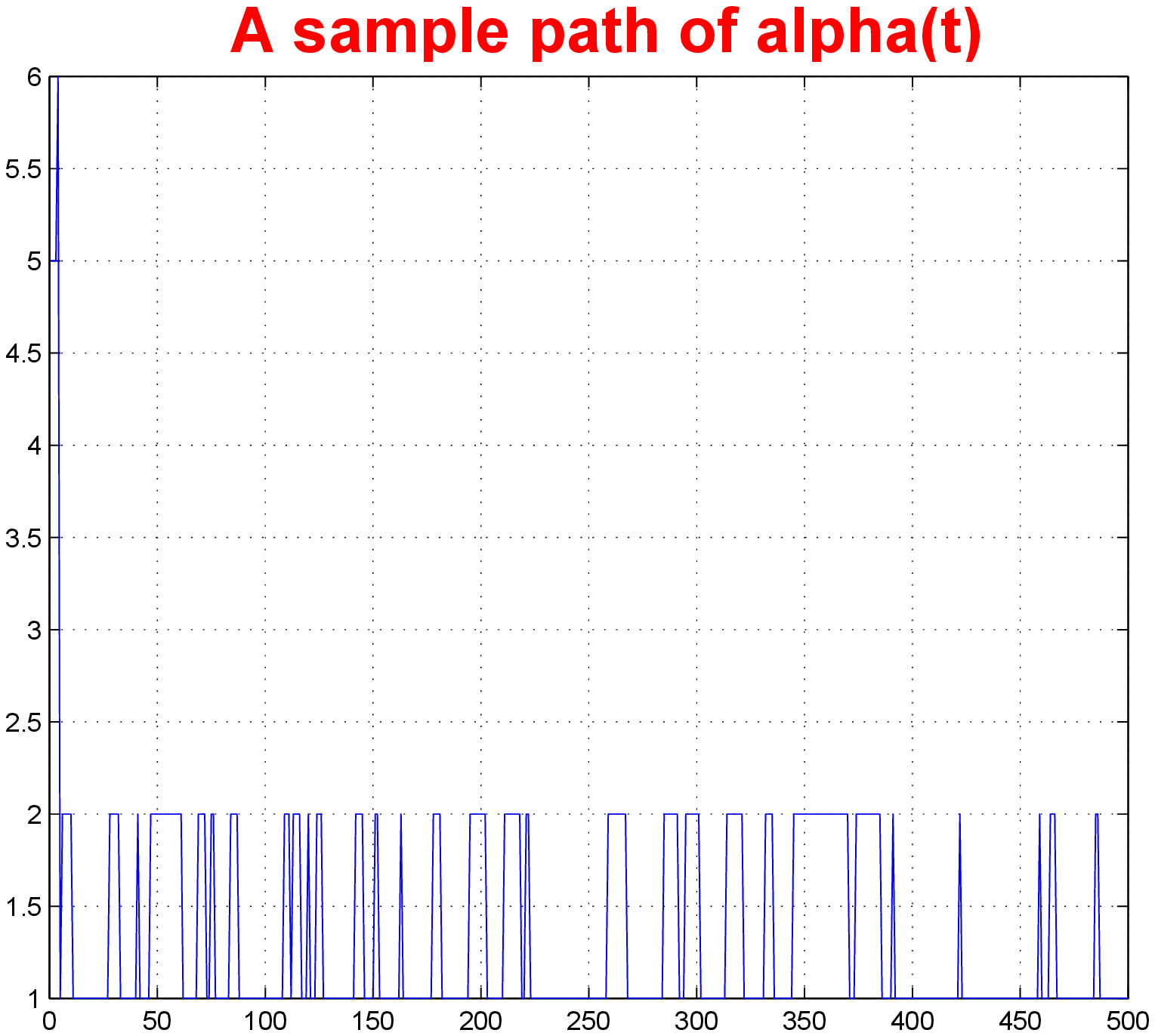}
\includegraphics[width=80mm]{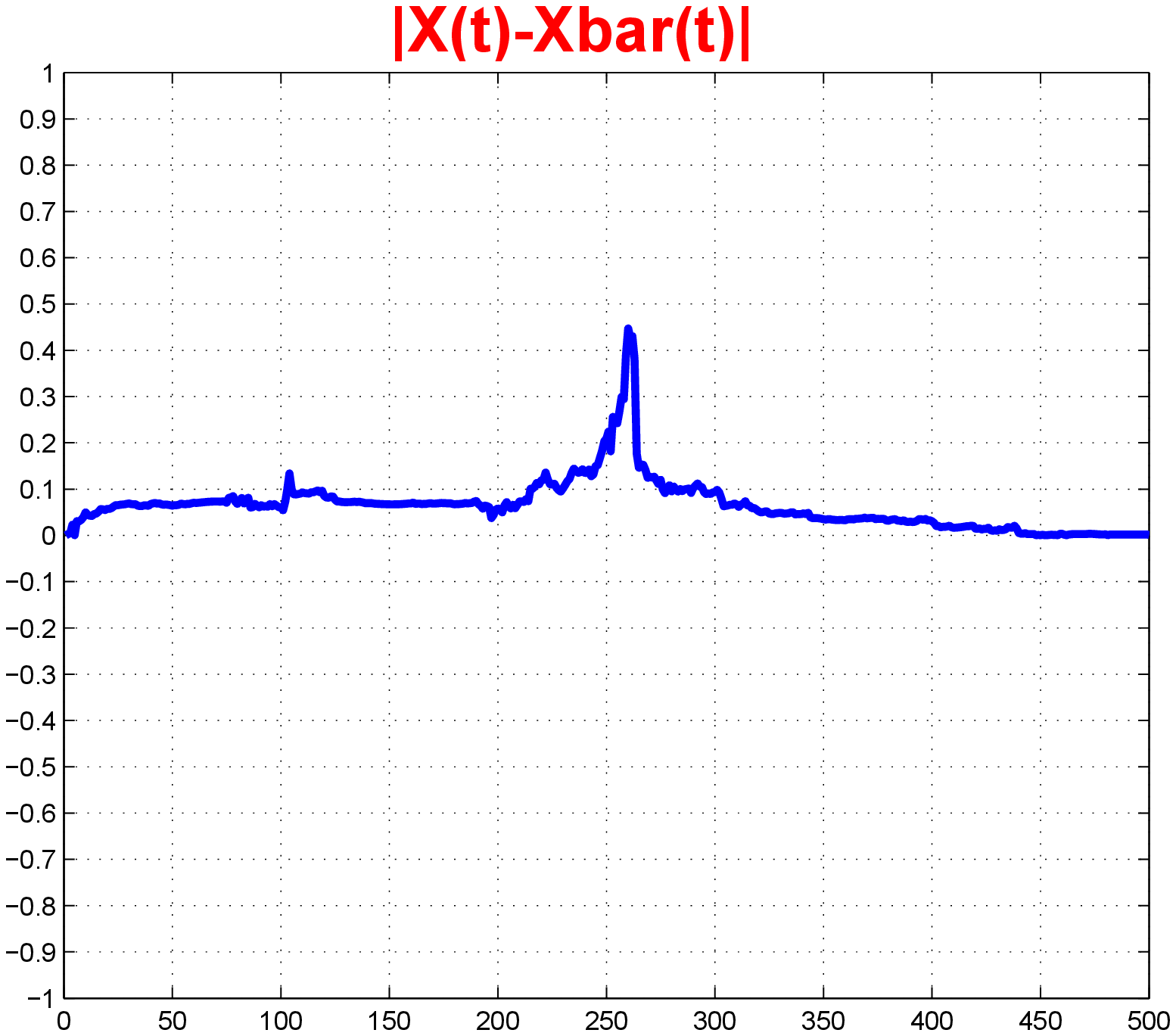}
\caption{Sample path with $\e=0.1$ in \exmref{6.2}}
\end{figure}
\begin{figure}
\includegraphics[width=80mm]{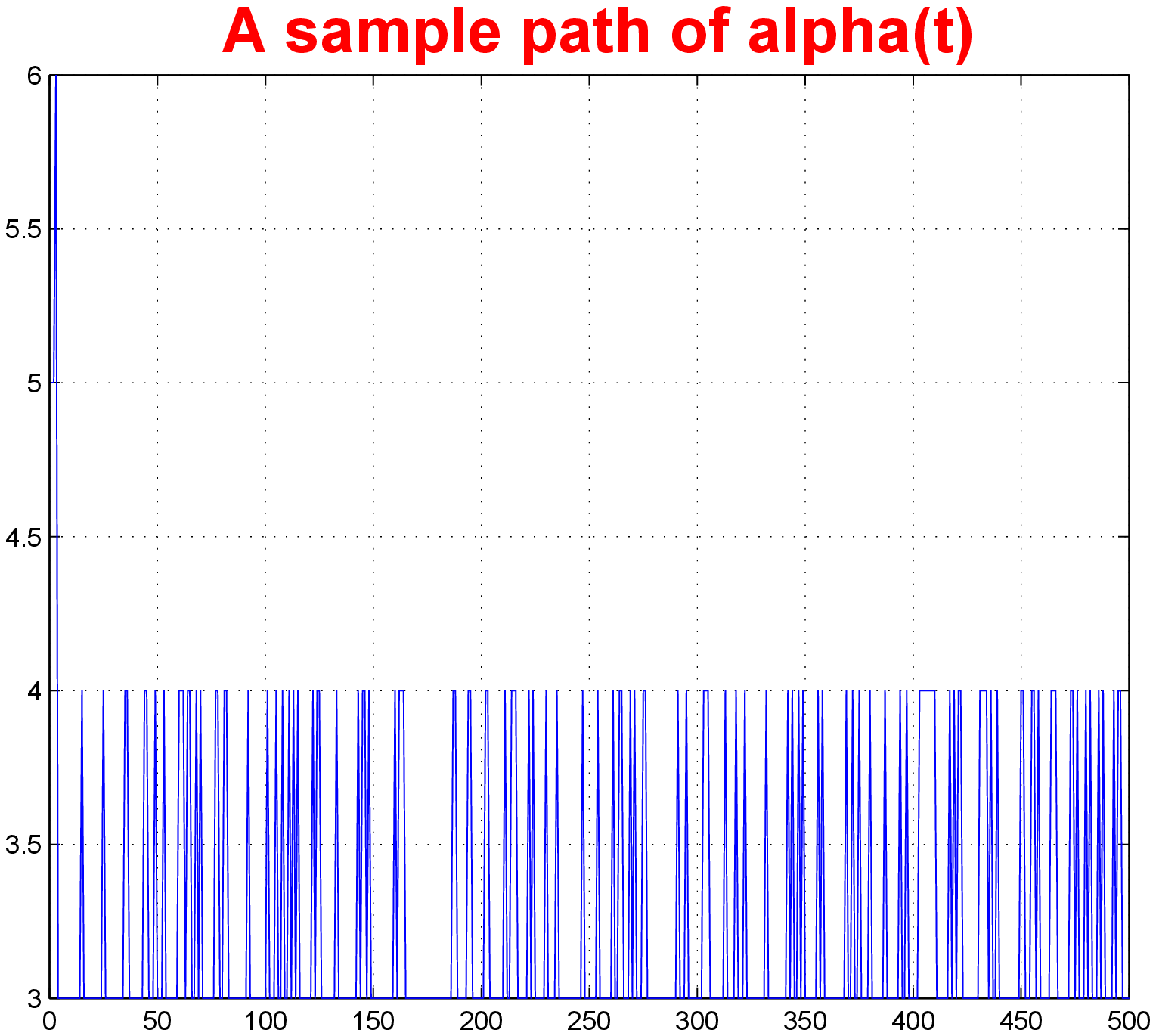}
\includegraphics[width=80mm]{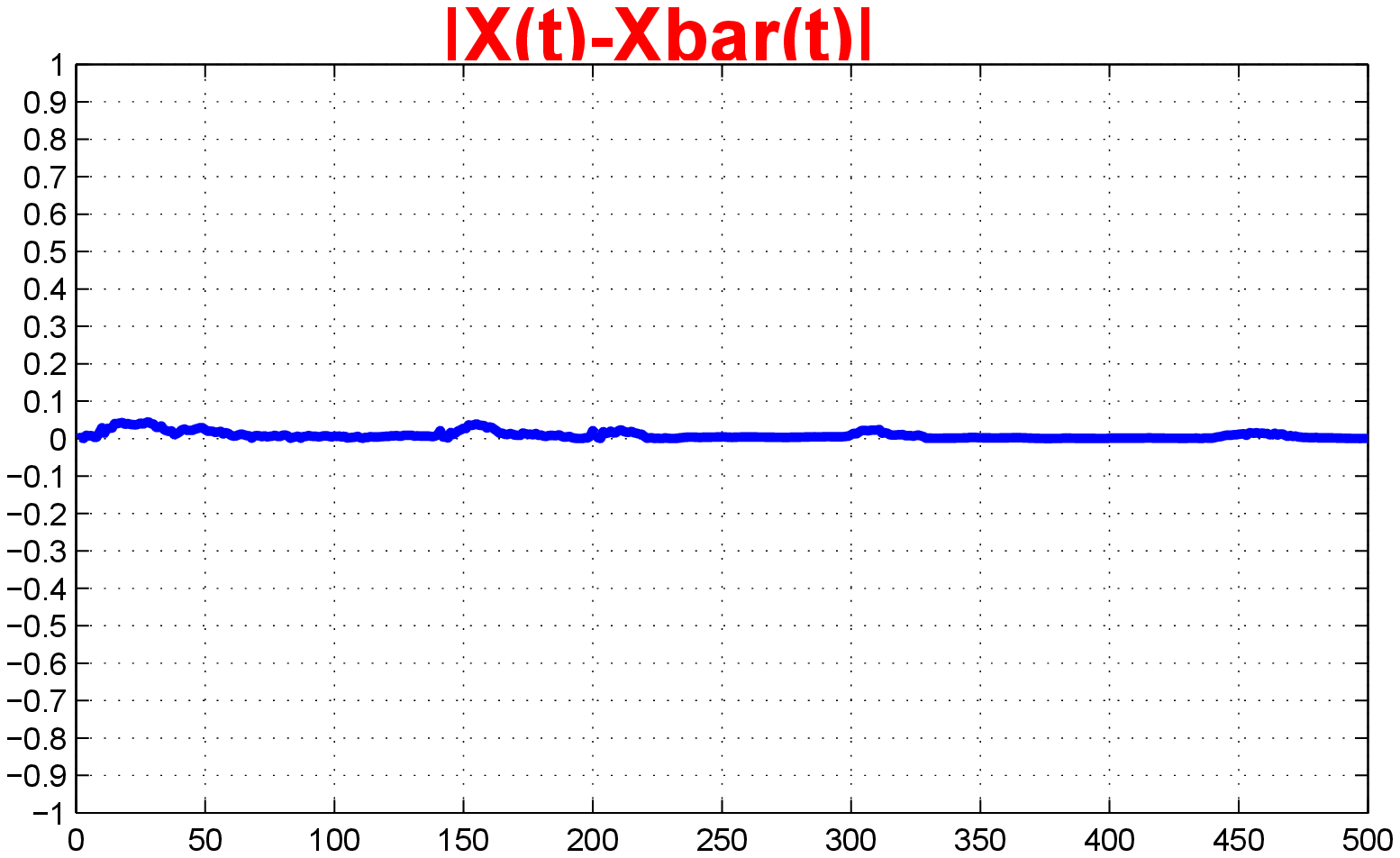}
\caption{Sample path with $\e=0.01$ in \exmref{6.2}}
\end{figure}
\begin{exm}\label{6.2} {\rm  We consider the Markov chain
$\al^\e(t)\in\mathcal{M}=\{1,2,3,4,5,6\}$ whose generator
$Q^\e$ is given by \eqref{mar} with
\bea \ad \wdt Q =\left(
\begin{array}{rrrrrr}
  -1 & 1 & 0 & 0 & 0 & 0 \\
   2 & -2 & 0 & 0 & 0 & 0 \\
 0 & 0 & -1 & 1 & 0 & 0 \\
 0 & 0 & 3 &-3 & 0 & 0 \\
 0 & 0 & 1 & 1 & -3 & 1 \\
 1 & 0 & 0 & 1 & 1 &-3 \\
 \end{array} \right), \\
 \ad \wdh Q=  \left( \begin{array}{rrrrrr}
  -1 & 1 & 0 & 0 & 0 & 0\\
 2 & -2 & 0 & 0 & 0 & 0\\
 0& 0 &-1 & 1 &0&0\\
 0 & 0 & 3 & -3 & 0& 0\\
 0 & 0 & 1 & 1 & -3 & 1\\
 1 & 0 &0 & 1 &1 &-3\\
 \end{array} \right)
.\eea
We use $x^\e(0)=0, r(t,1)=.5, r(t,2)=-.1, r(t,3)=.5, r(t,4)=-.1,
r(t,5)=.2, r(t,6)=.4, B(t,1)=1, B(t,2)=2, B(t,3)=-1, B(t,4)=-2,
B(t,5)=1, B(t,6)=2,
\sigma(t,1)=\sigma(t,2)=\sigma(t,3)=\sigma(t,4)=\sigma(t,5)=\sigma(t,6)=1$.
Sample paths of $\al^\e(t)$, trajectories of $|x^\e(t)-\lbar{x}(t)|$
are given in Figure 3 for $\e=.1$ and in Figure 4 for $\e=.01$. We
omit the error bounds here yet the result is similar to
\exmref{6.1}. }\end{exm}

It can be seen from the two graphs that the smaller the $\e$ the
more rapidly $\al^\e(\cdot)$ jumps, which results in better
approximations.

 \section{Further Remarks}\label{sec:fur}
 Motivated by platoon control systems,
 this work establishes
 a weak convergence property that  leads to a limit problem of much reduced complexity for
 mean-variance type of control under randomly regime switching systems. Our methodology uses
 a two-time-scale formulation to relate the underlying problem
 with that of the limit problem.
 Accompanying our recent work \cite{YYWZ}, this paper also demonstrates the
 near-optimal controls using numerical examples.
  Our approach provides a
 systematic method to reduce the complexity of the underlying
 system. In lieu of handling large dimensional systems, we
  need only solve a reduced set of limit equations that have much
 smaller dimensions.
 Future research efforts can be directed to the study of non-definite control problems in the hybrid systems,
 in which the Markov chain is a hidden process.
 Then a
 Wonham filter may be developed.
Another direction is to look into the possibility of treating
distributed controls with built-in communication complexity measures.
All of these deserve
 more thoughts and further considerations.

\end{document}